\documentclass[11pt]{article}

\usepackage{url,epsf,amsmath,amssymb,amsfonts,graphicx,tikz}
\usepackage{color}
\usepackage{hyperref}

\newtheorem{theorem}{\bf Theorem}[section]
\newtheorem{corollary}[theorem]{\bf Corollary}

\newtheorem{proposition}[theorem]{\bf Proposition}

\newtheorem{remark}[theorem]{\bf Remark}
\newtheorem{problem}[theorem]{\bf Problem}

\newcommand{\proof}{\noindent{\bf Proof.\ }}
\newcommand{\qed}{\hfill $\Box$ \bigskip}

\newcommand{\sr}{{}_{\rm SR}}
\newcommand{\gp}{{\rm gp}}

\newcommand{\diam}{{\rm diam}}

\begin{document}

\title{The general position problem and strong resolving graphs}

\author{
Sandi Klav\v zar $^{a,b,c}$
\and
Ismael G. Yero $^d$
}

\date{}

\maketitle

\begin{center}
$^a$ Faculty of Mathematics and Physics, University of Ljubljana, Slovenia \\
\medskip

$^b$ Faculty of Natural Sciences and Mathematics, University of Maribor, Slovenia \\
\medskip

$^{c}$ Institute of Mathematics, Physics and Mechanics, Ljubljana, Slovenia \\
\medskip

$^{d}$ Departamento de Matem\'aticas, Universidad de C\'adiz, Algeciras, Spain
\end{center}

\begin{abstract}
The general position number ${\rm gp}(G)$ of a connected graph $G$ is the cardinality of a largest set $S$ of vertices such that no three pairwise distinct vertices from $S$ lie on a common geodesic. It is proved that ${\rm gp}(G)\ge \omega(G_{\rm SR}$, where $G_{\rm SR}$ is the strong resolving graph of $G$, and $\omega(G_{\rm SR})$ is its clique number. That the bound is sharp is demonstrated with numerous constructions including for instance direct products of complete graphs and different  families of strong products, of generalized lexicographic products, and of rooted product graphs. For the strong product it is proved that $\gp(G\boxtimes H) \ge \gp(G)\gp(H)$, and asked whether the equality holds for arbitrary connected graphs $G$ and $H$. It is proved that the answer is in particular positive for strong products with a complete factor, for strong products of complete bipartite graphs, and for certain strong cylinders.
\end{abstract}

\noindent
{\bf E-mails}: sandi.klavzar@fmf.uni-lj.si, ismael.gonzalez@uca.es


\medskip\noindent
{\bf AMS Subj. Class.}: 05C12, 05C76

\section{Introduction}
\label{sec:intro}

The general position problem was recently and independently introduced in~\cite{manuel-2018a, ullas-2016}. If $G = (V(G), E(G))$ is a graph, then $S\subseteq V(G)$ is a {\em general position set} if no triple of vertices from $S$ lie on a common geodesic in $G$. The {\em general position problem} is to find a largest general position set of $G$, the order of such a set is the {\em general position number} $\gp(G)$ of $G$. A general position set of $G$ of order $\gp(G)$ is shortly called {\em gp-set}. The general position problem has been further studied in a sequence of very recent papers~\cite{anand-2019+, ghorbani-2019, manuel-2018b, patkos-2019+}.

A vertex $u$ of a connected graph $G$ is \emph{maximally distant} from a vertex $v$ if every $w\in N(u)$ satisfies $d_G(v,w)\le d_G(u,v)$, where $N(u)$ is the open neighborhood of $u$. If $u$ is maximally distant from $v$, and $v$ is maximally distant from $u$, then $u$ and $v$ are \emph{mutually maximally distant} (MMD for short). The {\em strong resolving graph} $G\sr$ of $G$ has $V(G)$ as the vertex set, two vertices being adjacent in $G\sr$ if they are MMD in $G$. The notion of the strong resolving graph was introduced in~\cite{Oellermann2007} as a tool to study the strong metric dimension. It was proved that the problem of finding the strong metric dimension of a graph $G$ can be transformed to the problem of finding the vertex cover number of  $G\sr$. Further on, the strong resolving graph itself was remarked as a kind of graph transformation in~\cite{Kuziak-2018}, and several characterizations and realizations of it were described.

Now, one of the open problems presented in~\cite{Kuziak-2018} concerns finding applications for the strong resolving graph construction, other than that of computing the strong metric dimension of graphs. In this paper we give a partial answer to this problem by establishing a connection between the general position number of a graph $G$ and the clique number of the graph $G\sr$. More precisely, in  Theorem~\ref{thm:lower-bound} we prove that $\gp(G)\ge \omega(G\sr)$ holds for any connected graph $G$. Then we demonstrate with different infinite families of graphs, including direct products of complete graphs, that the bound is sharp. We also show that for any integers $r\ge t\ge 2$, there exists a graph $G$ such that $\gp(G)=r$ and $\omega(G\sr)=t$. In Section~\ref{sec:strong} we focus on strong products of graphs. We prove that $\gp(G\boxtimes H) \ge \gp(G)\gp(H)$ holds for connected graphs $G$ and $H$ and that the bound is again sharp. In particular, if $\gp(G) = \omega(G\sr)$, then $\gp(G\boxtimes K_n) = n\cdot \gp(G) = \omega((G\boxtimes K_n)\sr)$. We close the section with a question on whether actually the equality $\gp(G\boxtimes H) = \gp(G)\gp(H)$ holds for arbitrary connected graphs $G$ and $H$. In Section~\ref{sec:lexico} we give additional large families of graphs, based on the generalized lexicographic product, for which the equality in  Theorem~\ref{thm:lower-bound} holds. In the final section we determine the general position number for different rooted product graphs and relate the values with the corresponding clique numbers of strong resolving graphs.

Before giving our results, we list in the next section definitions and concepts not yet given, as well as some results needed later.

\section{Preliminaries}
\label{sec:prelim}

For a positive integer $k$ we will use the notation $[k] = \{1,\ldots, k\}$. If $G = (V(G), E(G))$ is a graph, then $n(G) = |V(G)|$ and $m(G) = |E(G)|$. If $X\subseteq V(G)$, the the subgraph of $G$ induced by $X$ is denoted $\langle X\rangle$.

The distance $d_G(u,v)$ between vertices $u$ and $v$ of a graph $G$ is the number of edges on a shortest $u,v$-path. A subgraph $H$ of a graph $G$ is {\em isometric} if $d_H(u,v) = d_G(u,v)$ holds for all $u,v\in V(H)$. A set of subgraphs $\{H_1,\ldots, H_k\}$ of a graph $G$ is an {\em isometric cover} of $G$ if each $H_i$, $i\in [k]$, is isometric in $G$ and $\bigcup_{i=1}^k V(H_i) = V(G)$. With this concept in hand we can recall the following.

\begin{theorem} {\rm \cite[Theorem 3.1]{manuel-2018a}}
\label{thm:isometric-cover}
If $\{H_1,\ldots, H_k\}$ is an isometric cover of $G$, then
$$\gp(G) \le \sum_{i=1}^k \gp(H_i)\,.$$
\end{theorem}

If $G$ is a connected graph, $S\subseteq V(G)$, and ${\cal P} = \{S_1, \ldots, S_p\}$ a partition of $S$, then ${\cal P}$ is \emph{distance-constant} (alias ``distance-regular''~ \cite[p.~331]{kante-2017}) if for any $i,j\in [p]$, $i\ne j$, the distance $d_G(u,v)$, where $u\in S_i$ and $v\in S_j$, is independent of the selection of $u$ and $v$. This distance is then the distance $d_G(S_i,S_j)$ between the parts $S_i$ and $S_j$. A distance-constant partition ${\cal P}$ is {\em in-transitive} if $d_G(S_i, S_k) \ne d_G(S_i, S_j) + d_G(S_j,S_k)$ holds for pairwise different indices $i,j,k\in [p]$. With these concepts, general position sets can be characterized as follows.

\begin{theorem} {\rm \cite[Theorem 3.1]{anand-2019+}}
\label{thm:gpsets}
Let $G$ be a connected graph. Then $S\subseteq V(G)$ is a general position set if and only if the components of $\langle S\rangle$ are complete subgraphs, the vertices of which form an in-transitive, distance-constant partition of $S$.
\end{theorem}

Let $\eta(G)$ denote the maximum order of an induced complete multipartite subgraph of the complement $\overline{G}$ of a graph $G$. Then we have:

\begin{theorem} {\rm \cite[Theorem 4.1]{anand-2019+}}
\label{thm:diameter2}
If $\diam(G) = 2$, then $\gp(G) = \max\{\omega(G), \eta(G)\}$.
\end{theorem}

Vertices $u$ and $v$ of a graph $G$ are {\em true twins} if $N[u] = N[v]$, where $N[u]$ is the closed neighborhood of $u$. Note that only adjacent vertices can be true twins.

\begin{proposition}
\label{prop:diam2-twin-free}
If $G$ has no true twins and $\diam(G) = 2$, then $\gp(G) = \omega(G\sr)$ if and only if $\gp(G) = \alpha(G)$.
\end{proposition}

\proof
Since $G$ is true-twin free, $G\sr$ is isomorphic to the complement $\overline{G}$ of $G$. Hence $\omega(G\sr) = \alpha(G)$ and thus the conclusion.
\qed

The Petersen graph $P$ is a sporadic example of a graph without tree twins and of diameter $2$ for which $\gp(G) \ne \omega(G\sr)$. Indeed, $\gp(P) = 6$ and $\omega(P\sr) = \alpha(P) = 4$.

Let $G$ and $H$ be graphs. Among the four standard graph products, we will consider the {\em direct product} $G\times H$, the {\em strong product} $G\boxtimes H$, and  the {\em lexicographic product} $G[H]$. The vertex set of all these products is $V(G)\times V(H)$. Let $(g,h), (g',h')\in V(G)\times V(H)$. In $G\times H$, the vertices $(g,h)$ and $(g',h')$ are adjacent if $gg\in E(G)$ and $hh'\in E(H)$. In $G\boxtimes H$, the vertices $(g,h)$ and $(g',h')$ are adjacent if one of the following three conditions hold: (i) $gg\in E(G)$ and $h=h'$,  (ii) $g=g$ and $hh'\in E(H)$, (iii) $gg\in E(G)$ and $hh'\in E(H)$. Finally, in $G[H]$ the vertices $(g,h)$ and $(g',h')$ are adjacent if either $gg\in E(G)$, or $g=g'$ and $hh'\in E(H)$. We note that the lexicographic product is also denoted with $G\circ H$ to emphasize the associativity of the operation, but here we use $G[H]$ to be consistent with the generalized lexicographic product (to be defined below). If $G\ast H$ is one of the above products and $h\in V(H)$, then the subgraph of $G\ast H$ induced by $\{(g,h):\ g\in V(G)\}$ is called a {\em $G$-layer}. Analogously $H$-layers are defined. In $G\times H$, each $G$-layer is an edgeless graph of order $n(G)$. In all other above products,
each $G$-layer is isomorphic to $G$. For more information on the standard graph products see the book~\cite{hik-2011}, here we just recall the following well-known result (cf.~\cite[Proposition 5.4]{hik-2011}) needed later.

\begin{proposition}
\label{prop:distance-in-strong}
If $(g,h)$ and $(g',h')$ are vertices of a strong product $G \boxtimes H$, then
$$d_{G \boxtimes H}\left((g,h),(g',h')\right) = \max\{d_G(g,g'), d_H(h,h')\}\,.$$
\end{proposition}

\section{The lower bound and equality cases}
\label{sec:lower-bound}

In this section we first prove the key result that connects the general position problem with the strong resolving graphs.

\begin{theorem}
\label{thm:lower-bound}
If $G$ is a connected graph, then $\gp(G)\ge \omega(G\sr)$. Moreover, equality holds if and only if $G$ contains a gp-set that induces a complete subgraph of $G\sr$.
\end{theorem}

\proof
Let $S\subseteq V(G\sr)$ induce a complete subgraph of $G\sr$. This means that any two vertices $x,y\in S$ are MMD in $G$. We now consider the vertices of $S$ in the graph $G$. If there are three distinct vertices $x,y,z\in S$ lying on a common geodesic, say $y$ lies in an $x,z$-geodesic, then neither $x,y$ nor $y,z$ are MMD in $G$, which is a contradiction. Thus, any three vertices of $S$ do not lie in a common geodesic of $G$, and therefore, $S$ is a general position set in $G$. Selecting $S$ to be a complete subgraph $G\sr$ of order $\omega(G\sr)$ leads to the desired bound.

Suppose now that $\gp(G) = \omega(G\sr)$. By the above, any complete subgraph of $G\sr$ of order $\omega(G\sr)$ yields a gp-set. Conversely, let $S$ be a gp-set of $G$ that forms a complete subgraph of $G\sr$. Then, using the already proven inequality $\gp(G)\ge \omega(G\sr)$, we have
$$|S| \le \omega(G\sr) \le \gp(G) = |S|\,,$$
from which we conclude that $\omega(G\sr) = \gp(G)$.
\qed

One would immediately think of characterizing the class of graphs achieving the equality in Theorem~\ref{thm:lower-bound}. However, such a characterization seems to be elusive because of the great variety of different structures that can appear. In the following we justify this variety and begin a couple of simple examples that were implicitly known previously.

\begin{itemize}
\item Block graphs, in particular complete graphs and trees. \\
Indeed, in~\cite{manuel-2018a} it was observed that in block graphs the set of simplicial vertices forms a gp-set. Since simplicial vertices of a graph $G$ also form a set of MMD vertices of a graph (equivalently, they form a complete subgraph of $G\sr$), Theorem~\ref{thm:lower-bound} implies that $\gp(G) = \omega(G\sr)$ if $G$ is a block graph.

\item Complete multipartite graphs. \\
Let $G = K_{n_1,\ldots, n_k}$, where $k\ge 2$ and $n_1\ge n_2\ge \cdots \ge n_k\ge 2$. Then it is easy to see that the vertices of the $n_1$-partite set form a maximum general position set. Moreover, the vertices of this set also form a set of mutually maximally distant vertices of $G$. Hence, $\gp(G) = \omega(G\sr)$ by Theorem~\ref{thm:lower-bound}.
\end{itemize}

Let $G$ and $H$ be graphs where $V(G) = \{v_1, \ldots ,v_{n}\}$. The {\em corona} $G\odot H$ of graphs $G$ and $H$ is obtained from the disjoint union of $G$ and $n$ disjoint copies of $H$, say $H_1,\ldots, H_{n}$, where for all $i\in [n]$, the vertex $v_i\in V(G)$ is adjacent to each vertex of $H_i$. Then we have another equality case:

\begin{proposition}
\label{prop:corona}
If $H=\bigcup_i K_{n_i}$, $n_i\ge 1$, then for every graph $G$, $\gp(G\odot H) = \omega((G\odot H)\sr)$.
\end{proposition}

\proof
From \cite[Theorem 4.3]{ghorbani-2019}, it can be noticed that $\gp(G\odot H) = n(G)\sum_i n_i$, and also that the union of the sets of vertices of every copy of $H$ in $G\odot H$ form a gp-set $S$ of $G\odot H$. Every two vertices belonging to one copy of $H$ are MMD, as well as are MMD every two vertices belonging to two different copies of $H$. Hence $S$ forms a complete subgraph of $(G\odot H)\sr$. Thus we deduce the equality by Theorem~\ref{thm:lower-bound}.
\qed

We note in passing that Proposition~\ref{prop:corona} remains valid in a more general setting when different disjoint unions of complete graphs are attached to the vertices of $G$.

In the next result we provide a family of direct product graphs for which the equality holds in Theorem~\ref{thm:lower-bound}.

\begin{proposition}
\label{prop:direct-complete}
If $n_1\ge n_2\ge 3$, then
$$\gp(K_{n_1}\times K_{n_2}) = \omega((K_{n_1}\times K_{n_2})\sr) = n_1 = \alpha((K_{n_1}\times K_{n_2})\sr)\,.$$
\end{proposition}

\proof
We first note that $\omega(K_{n_1}\times K_{n_2})=\min\{n_1,n_2\}=n_2$. On the other hand, since every two vertices of $K_{n_1}\times K_{n_2}$ belonging to two different copies of $K_{n_1}$ and of $K_{n_2}$ are adjacent, every maximal induced complete multipartite subgraph of $\overline{K_{n_1}\times K_{n_2}}$ is formed by the set of vertices of one copy of $K_{n_1}$ or of $K_{n_2}$. Thus, $\eta(K_{n_1}\times K_{n_2})=\max\{n_1,n_2\}=n_1$. Now, since $n_1\ge n_2\ge 3$, it follows that $\diam(K_{n_1}\times K_{n_2}) = 2$ and hence Theorem~\ref{thm:diameter2} yields $\gp(K_{n_1}\times K_{n_2})=\max\{\eta(K_{n_1}\times K_{n_2}),\omega(K_{n_1}\times K_{n_2})\}=n_1$. From~\cite{Kuziak-2018, Rodriguez-2014} it is known that $(K_{n_1}\times K_{n_2})\sr\cong K_{n_1}\Box K_{n_2}$ and since $\omega(K_{n_1}\Box K_{n_2})=\max\{n_1,n_2\}=n_1$, the first two equalities follows. The last equality then follows by Proposition~\ref{prop:diam2-twin-free}.
\qed

Note that if we consider $n_1>n_2=2$ in the result above, then $K_{n_1}\times K_2$ is of diameter $3$, and its strong resolving graph is $\overline{K}_{n_1}\Box K_{2}$. Thus, $\omega((K_{n_1}\times K_{2})\sr)=2$. Since $\gp(K_{n_1}\times K_{2})=n_1>2$, there is no equality as in the proposition.

Another example of direct products for which the equality in Theorem~\ref{thm:lower-bound} does not hold is $K_{r,t}\times K_n$, where $r\ge t\ge 2$ and $n\ge 3$. Since $\diam(K_{r,t}\times K_n) = 3$, \cite[Theorem 5.1]{anand-2019+} implies that $\gp(K_{r,t}\times K_n) = \alpha(K_{r,t}\times K_n)$. Since it is not difficult to verify that $\alpha(K_{r,t}\times K_n) = rn$, we get $\gp(K_{r,t}\times K_n) = rn$. On the other hand, from~\cite[Theorem 35]{Kuziak-2018} we know that $(K_{r,t}\times K_n)\sr\cong \bigcup_{i=1}^n K_{r+t}$, and so $\omega((K_{r,t}\times K_n)\sr)=r+t$. As $r\ge t\ge 2$ and $n\ge 3$ we have $rn > r + t$.

Based on the above special cases we pose the following question about a possible dichotomy in direct product.

\begin{problem}
Is it true that $\gp(G\times H) = \omega((G\times H)\sr)$ can only hold in the case when $\diam(G\times H) = 2$?
\end{problem}

To conclude the section we give the following realization result which intuitively indicates that one cannot expect some natural upper bound on $\gp(G)$ in terms of $\omega(G\sr)$.

\begin{proposition}
For any integers $r\ge t\ge 2$, there exists a graph $G$ such that $\gp(G)=r$ and $\omega(G\sr)=t$.
\end{proposition}

\proof
Since $r\ge t$, there exists a non-negative integer $q$ such that $r=t+q$. We now consider a graph $G_q$ defined as follows. We begin with $q$ copies of the cycle graph $C_4$ and $t-q$ copies of the graph $P_2$. Then we add an extra vertex $z$ and one edge between $z$ and exactly one vertex of each copy of $C_4$ and of $P_2$. We observe that the components of the strong resolving graph $(G_q)\sr$ are: one  complete graph of order $t$, $q$ complete graphs $K_2$, and $t+1$ isolated vertices. Thus $\omega((G_q)\sr)=t$. On the other hand, a set formed by two non-adjacent vertices of each copy of the cycle $C_4$ (those ones not adjacent to $z$), and one vertex of each copy of the path $P_2$, used to construct $G_q$, is a general position set of $G_q$, and so, $\gp(G_q)\ge 2q+t-q=t+q=r$. We can readily observe that such set is indeed a gp-set of $G_q$, and therefore $\gp(G)=r$, which completes the proof.
\qed

\section{Strong products}
\label{sec:strong}

If $G$ and $H$ are connected graphs, then each $G$-layer and each $H$-layer of $G\boxtimes H$ is an isometric subgraph of $G\boxtimes H$. Hence Theorem~\ref{thm:isometric-cover} gives the following upper bound.

\begin{corollary}
\label{cor:strong-upper}
If $G$ and $H$ are connected graphs, then
$$\gp(G\boxtimes H) \le \min\{n(G)\gp(H), n(H)\gp(G)\}\,.$$
\end{corollary}

We will later see that the bound of Corollary~\ref{cor:strong-upper} is tight. On the other hand we have the following lower bound.

\begin{theorem}
\label{thm:strong-lower}
If $G$ and $H$ are connected graphs, then $\gp(G\boxtimes H) \ge \gp(G)\gp(H)$.
\end{theorem}

\proof
Let $S_G$ and $S_H$ be gp-sets of $G$ and $H$, respectively, so that $|S_G| = \gp(G)$ and $|S_H| = \gp(H)$. We claim that $S_G\times S_H$ is a general position set of $G\boxtimes H$. To prove it, consider arbitrary pairwise different vertices of $S_G\times S_H$, say $(g, h)$, $(g', h')$, $(g'', h'')$, and assume on the contrary that in $G\boxtimes H$ there exists a shortest $(g,h),(g'',h'')$-path $P$ that passes through $(g',h')$. We now distinguish several cases.

Suppose first that $g=g'=g''$. Since $(g, h)$, $(g', h')$, $(g'', h'')$ are pairwise different vertices of $G\boxtimes H$, the vertices $h$, $h'$, $h''$ are then pairwise different. But then the projection of $P$ to $H$ is a shortest $h,h''$-path that contains $h'$, a contradiction. Similarly, if $g$, $g'$, $g''$ are pairwise different, then the projection of $P$ to $G$ is a shortest $g,g''$-path that contains $g'$.

Suppose next that $g=g'$ and $g''\ne g$. Then clearly $h\ne h'$. If $h''$ is different from both $h$ and $h'$, then, as above, consider the projection of $P$ to $H$ to get a contradiction. The other subcase is that $h'' = h'$ (the subcase $h'' = h$ is treated analogously). Let $d_G(g,g'') = k$ and $d_H(h,h'') = \ell$. By Proposition~\ref{prop:distance-in-strong}, $d_{G\boxtimes H}((g,h), (g'',h'')) = \max\{k,\ell\}$. Denoting by $P'$ the $(g,h),(g,h')$-subpath of $P$ and by $P''$ the $(g,h'),(g'',h'')$-subpath of $P$, we get that $\max\{k,\ell\} = |P| = |P'| + |P''| \ge \ell + k$, a contradiction since $k\ge 1$ and $\ell \ge 1$. The case $g =  g''$, $g' \ne g$, and the case $g'=g''$, $g\ne g'$, are treated analogously.
\qed

In~\cite[Theorem 3.3]{manuel-2018b} it was proved that $\gp(P_\infty \boxtimes P_\infty) = 4$. Since the strong grid $P_n\boxtimes P_m$ is an isometric subgraph of $P_\infty \boxtimes P_\infty$ for each $n,m\ge 2$, it follows that $\gp(P_n \boxtimes P_m) \le 4$. On the other hand, as $P_n \boxtimes P_m$ contains $K_4$ we also have $\gp(P_n \boxtimes P_m) \ge 4$. We conclude that
\begin{equation}
\label{eq:Pn strong Pm}
\gp(P_n \boxtimes P_m) = 4,\  n,m \ge 2\,.
\end{equation}
This result shows that the bound in Theorem~\ref{thm:lower-bound} is sharp. More sharpness examples are provided with the next result which also shows the tightness of Corollary~\ref{cor:strong-upper}.

\begin{proposition}
\label{prop:strong-one-factor-complete}
If $G$ is a connected graph and $n\ge 1$, then $\gp(G\boxtimes K_n) = n\cdot \gp(G)$. Moreover, if $\gp(G) = \omega(G\sr)$, then $\gp(G\boxtimes K_n) = \omega((G\boxtimes K_n)\sr)$.
\end{proposition}

\proof
The first assertion follows by combining Theorem~\ref{thm:strong-lower} with Corollary~\ref{cor:strong-upper}.

Suppose now that in addition $\gp(G) = \omega(G\sr)$ holds. In view of Theorem~\ref{thm:lower-bound}, there exists a gp-set $S_G$ of $G$ that induces a complete subgraph of $G\sr$. By the proof of Theorem~\ref{thm:strong-lower}, $S_G\times V(K_n)$ is a gp-set of $G\boxtimes K_n$. The components of the subgraph of $G\boxtimes K_n$ induced by $S_G\times V(K_n)$ are of the form $Q\boxtimes K_n$, where $Q$ is a complete component induced by $S_G$. If $(g,x)$ and $(g',x')$ belong to different components $Q\boxtimes K_n$ and $Q'\boxtimes K_n$ induced by $S_G\times V(K_n)$, then with Proposition~\ref{prop:distance-in-strong} in mind we have $d_{G\boxtimes K_n}((g,x), (g',x')) = d_{G}(g, g)$. Since $g$ and $g'$ are MMD, this implies that also $(g,x)$ and $(g',x')$ are MMD. Suppose next that $(g,x)$ and $(g',x')$ belong to the same component $Q\boxtimes K_n$. If $g = g'$, then $(g,x)$ and $(g',x')$ are clearly MMD. Suppose now that $g\ne g'$. Then, since $g$ and $g'$ are adjacent and MMD in $G$, the vertices $g$ and $g'$ are true twins. But then it follows that $(g,x)$ and $(g',x')$ are MMD in $G\boxtimes K_n$. The second assertion now follows from Theorem~\ref{thm:gpsets}.
\qed

Since $\gp(T)=t$ for every tree $T$ with $t$ leaves, we have $\gp(T\boxtimes P_n)\ge 2t$ by Theorem~\ref{thm:strong-lower}. We next show that this becomes an equality for an infinite number of trees. To this end, we say that a tree $T$ belongs to a family $\mathcal{T}$ if there exits a finite sequence $T_1,\ldots,T_r$, $r\ge 1$, of trees such that,
\begin{itemize}
  \item $T_1$ is a path on at least three vertices;
  \item $T_2$ is obtained from $T_1$ by adding a path $P$ of order at least $3$ and joining by an edge one not leaf vertex of $P$ with one not leaf vertex of $T_1$;
  \item for every $i\in \{3,\ldots,r\}$, $T_i$ is obtained from $T_{i-1}$ by adding a path $P$ of order at least $3$ and joining by an edge one not leaf vertex of $P$ with one vertex of degree larger than two of $T_{i-1}$; and
  \item $T = T_r$.
\end{itemize}
Note that if $T\in \mathcal{T}$ is obtained by the above construction in $r$ steps, then $T$ has exactly $2r$ leaves.

\begin{proposition}
\label{prop:special-trees}
If $T\in \mathcal{T}$ and has $2r$ leaves, then $\gp(T\boxtimes P_n) = 4r = \omega((T\boxtimes P_n)\sr)$.
\end{proposition}

\proof
From Theorem~\ref{thm:strong-lower} we get $\gp(T\boxtimes P_n)\ge 4r$. We next show that this is also the exact value.

Let $P_{n_i}$, $i\in [r]$, be the path used to generate $T$ in the $i^{\rm th}$ step of the construction of $T$. Let $S_i = V(P_{n_i}) \times V(P_n)$, $i\in [r]$. Then note that $S_1,\dots,S_r$ form a partition of $V(T\boxtimes P_n)$, where each $S_i$ induces a graph isomorphic to the strong grid graph $P_{n_i}\boxtimes P_n$.

Since $\{S_1,\dots,S_r\}$ form an isometric cover of $T\boxtimes P_n$, Theorem~\ref{thm:isometric-cover} and~\eqref{eq:Pn strong Pm} imply that
$$\gp(T\boxtimes P_n) \le \sum_{i=1}^r \gp(\langle S_i\rangle) = 4r\,,$$
hence the first equality follows.

From~\cite[Theorem 40]{Kuziak-2018} we know that $T\sr \boxtimes (P_n)\sr$ is a subgraph of $(T\boxtimes P_n)\sr$. Since $T\sr \boxtimes (P_n)\sr$ contains a clique of size $4r$, we then also have such a clique in $(T\boxtimes P_n)\sr$ and so $\omega((T\boxtimes P_n)\sr) \ge 4r$. Theorem~\ref{thm:lower-bound} completes the argument.
\qed

\begin{proposition}
\label{prop:strong-two-bipar}
If $r_1\ge t_1\ge 1$ and $r_2\ge t_2\ge 1$, then
$$\gp(K_{r_1,t_1}\boxtimes K_{r_2,t_2})=r_1r_2=\omega((K_{r_1,t_1}\boxtimes K_{r_2,t_2})\sr) = \alpha(K_{r_1,t_1}\boxtimes K_{r_2,t_2})\,.$$
\end{proposition}

\proof
We first observe that $\diam(K_{r_1,t_1}\boxtimes K_{r_2,t_2}) = 2$ and thus Theorem~\ref{thm:diameter2} applies.

The set obtained from the Cartesian product of the partite sets of cardinality $r_1$ and $r_2$ of $K_{r_1,t_1}$ and $K_{r_2,t_2}$, respectively, forms a maximal induced complete multipartite subgraph of $K_{r_1,t_1}\boxtimes K_{r_2,t_2}$ or cardinality $r_1r_2$. Since $\omega(K_{r_1,t_1}\boxtimes K_{r_2,t_2})=4$, we deduce that $\gp(K_{r_1,t_1}\boxtimes K_{r_2,t_2})=r_1r_2$.

On the other hand, since $K_{r_1,t_1}\boxtimes K_{r_2,t_2}$ has diameter two and has not true twin vertices, the strong resolving graph $(K_{r_1,t_1}\boxtimes K_{r_2,t_2})\sr$ is just the complement of $K_{r_1,t_1}\boxtimes K_{r_2,t_2}$. Thus, we obtain that $\omega((K_{r_1,t_1}\boxtimes K_{r_2,t_2})\sr)=\alpha(K_{r_1,t_1}\boxtimes K_{r_2,t_2})=r_1r_2$, hence the first two equalities.

The last equality follows by Proposition~\ref{prop:diam2-twin-free}.
\qed

\begin{theorem}
\label{th:odd-cylinder}
If $r\ge 2$ and $t\ge 1$, then $6\le \gp(P_r\boxtimes C_{2t+1})\le 7$. Moreover, if $t\in [2]$ or $r=2$, then $\gp(P_r\boxtimes C_{2t+1})=6$.
\end{theorem}

\proof
If $t=1$, then by Proposition~\ref{prop:strong-one-factor-complete}, $\gp(P_r\boxtimes C_{3})=\gp(P_r\boxtimes K_{3})= 3\gp(P_r)=6$. Hence, from now on we may assume $t\ge 2$.

Let $U=\{u_1,\dots,u_r\}$ and $V=\{v_1,\dots,v_{2t+1}\}$ be the vertex sets of $P_r$ and $C_{2t+1}$, respectively, with natural adjacencies. From Theorem \ref{thm:strong-lower}, we know that $\gp(P_r\boxtimes C_{2t+1})\ge 6$. A subpath $P$ of $C_{2t+1}$ which is of length at most $t$ is an isometric subgraph of $C_{2t+1}$, hence $U\times P$ induces an isometric subgraph of $P_r\boxtimes C_{2t+1}$. In particular this implies that the set $\{U\times \{v_1, \ldots, v_{t+1}\}, U\times \{v_{t+2}, \ldots, v_{2t+1}\}\}$ forms an isometric cover of $P_r\boxtimes C_{2t+1}$ consisting of two strong grids. Hence, again using Theorem~\ref{thm:isometric-cover} together with~\eqref{eq:Pn strong Pm} we infer that $\gp(P_r\boxtimes C_{2t+1})\le 8$.

We now suppose that $\gp(P_r\boxtimes C_{2t+1})=8$ and let $S$ be a gp-set of $P_r\boxtimes C_{2t+1}$. Let $S'$ be the projection of $S$ onto $C_{2t+1}$ and consider the following situations.

\medskip\noindent
{\bf Case 1}: $|S'|=8$. \\
This means that for each $v_i\in S'$ we have $|(U\times \{v_i\})\cap S|=1$. Without loss of generality we can assume that $v_1\in S'$. Consider now a partition of $V$ given by the sets $V_1=\{v_1,\dots,v_{t+1}\}$ and $V_2=\{v_{t+2},\dots,v_{2t+1}\}$. As noted above, $U\times V_1$ and $U\times V_2$ induce strong grids that are isometric subgraphs of $P_r\boxtimes C_{2t+1}$. Thus, by~\eqref{eq:Pn strong Pm} and since we have assumed $\gp(P_r\boxtimes C_{2t+1})=8$, we deduce $|S\cap (U\times V_1)|=4$ and $|S\cap (U\times V_2)|=4$. Analogously, if $V'_1=\{v_2,\dots,v_{t+1}\}$ and $V'_2=\{v_{t+2},\dots,v_{2t+1},v_1\}$, then also $U\times V'_1$ and $U\times V'_2$ induce two strong grids that are isometric subgraphs of $P_r\boxtimes C_{2t+1}$. Since $|S\cap (U\times V'_2)|=5$, we get a contradiction.

\medskip\noindent
{\bf Case 2}: $4\le |S'|\le 7$. \\
This means that there exists at least one vertex $v_i\in S'$ such that $|(U\times \{v_i\})\cap S|=2$. Note that for every $v_i\in S'$, it must happen $|(U\times \{v_i\})\cap S|\le 2$, otherwise we find a geodesic containing three vertices of $S$. Without loss of generality we can assume that $v_1\in S'$ satisfies that $|(U\times \{v_1\})\cap S|=2$. A similar argument as in Case 1 leads to a partition of $V$ given by $V''_1=\{v_2,\dots,v_{t+1}\}$ and $V''_2=\{v_{t+2},\dots,v_{2t+1},v_1\}$, and such that $U\times V''_1$ and $U\times V''_2$ induce strong grids that are isometric subgraphs of $P_r\boxtimes C_{2t+1}$ for which $|S\cap (U\times V'_2)|=6$, which is again not possible.

\medskip\noindent
{\bf Case 3}: $|S'|\le 3$. \\
Since for every $v_i\in S'$, it must happen $|(U\times \{v_i\})\cap S|\le 2$, we deduce that $|S|=\sum_{v_i\in S'}|(U\times \{v_i\})\cap S|\le 2|S'|\le 6$. This is a final contradiction proving that $|S| = 8$ is not possible.

We have thus proved that $\gp(P_r\boxtimes C_{2t+1})\le 7$. Let next $t = 2$. Then we consider again the projection $S'$ as defined above, but in this case we clearly have $|S'|\le 5$. Now, if $|S| \in \{7,8\}$, then we get a contradiction along the same lines as above.  Hence $\gp(P_r\boxtimes C_{5}) = 6$. Finally, if $r=2$, then the situation in which $|S'|=7$ leads to the existence of seven vertices lying in different layers of the factor graph $P_2$. But then there are three of such vertices lying on the same geodesic, which is not possible and so $\gp(P_2\boxtimes C_{2t+1})=6$.
\qed

Upper bounds on the general position number of the cylinder $P_r\boxtimes C_{2t}$  and of the torus $C_r\boxtimes C_t$, can be deduced by using similar techniques as in the proof above, except that in the last two cases we split the torus into two cylinders. On the other hand, lower bounds can be obtained from Theorem \ref{thm:strong-lower}. That is next stated.

\begin{remark}
\label{rem:bounds}
Let $r,t$ be two integers.
\begin{itemize}
  \item If $r\ge 2$ and $t\ge 3$, then $6\le \gp(P_r\boxtimes C_{2t})\le 8$.
  \item If $r\ge 5$ and $t\ge 3$, then $9\le \gp(C_r\boxtimes C_{2t})\le 16$.
  \item If $r\ge 4$ and $t\ge 2$, then $9\le \gp(C_r\boxtimes C_{2t+1})\le 14$.
\end{itemize}
\end{remark}

Using similar approach as in the proof of Theorem~\ref{th:odd-cylinder}, the last upper bound $14$ from Remark~\ref{rem:bounds} can be lowered to $13$.

Since $C_4$ is a complete bipartite graph and satisfies $\gp(C_4)=2$, from Proposition \ref{prop:strong-two-bipar}, we obtain that $\gp(C_4\boxtimes C_{4})=4$. Note that it also occurs the equality $\gp(C_4\boxtimes C_{4})=4=\omega((C_4\boxtimes C_{4})\sr)$ (for information on the structure of $(C_4\boxtimes C_{4})\sr$ see \cite{Kuziak-2018}).

Based on the results of this section we pose the following:

\begin{problem}
\label{prob:strong}
Is it true that if $G$ and $H$ are arbitrary connected graphs, then
$$\gp(G\boxtimes H) = \gp(G)\gp(H)\,?$$
\end{problem}

Assuming that the answer to the problem is positive, if $\gp(G) = \omega(G\sr)$ and $\gp(H) = \omega(H\sr)$, then $\gp(G\boxtimes H) = \omega((G\boxtimes H)\sr)$.

\section{Generalized lexicographic products}
\label{sec:lexico}

Let $G$ be a graph with $V(G) = \{g_1,\ldots, g_n\}$ and let $H_i$, $i\in [n]$, be pairwise disjoint graphs. Then the {\em generalized lexicographic product} $G[H_1,\ldots, H_n]$ has the vertex set
$$\bigcup_{i\in [n]}\{(g_i,h):\ h\in V(H_i)\}\,,$$
and the edge set
\begin{eqnarray*}
& & \{(g_i,h)(g_j,h'):\ g_ig_j\in E(G), h\in E(H_i), h'\in E(H_j)\} \cup \\
& & \bigcup_{i\in [n]} \{(g_i,h)(g_i,h'):\ hh'\in E(H_i)\}\,.
\end{eqnarray*}
In words, $G[H_1,\ldots, H_n]$ is obtained from $G$ by replacing each vertex $v_i\in V(G)$ with the graph $H_i$, and each edge $g_ig_{j}\in E(G)$ with all possible edges between $H_i$ and $H_{j}$. From this reason we will say that $v_i\in V(G)$ {\em expands} to $H_{i}$ in $G[H_1,\ldots, H_n]$.

The generalized lexicographic product was introduced by Sabidussi back in~\cite{sabidussi-1961}. If all the graphs $H_i$, $i\in [n]$, are isomorphic to a graph $H$, then the generalized lexicographic product $G[H_1,\ldots, H_n] = G[H,\ldots, H]$ becomes the standard lexicographic product $G[H]$.

\begin{theorem}
\label{thm:blow-up}
Let $G$ be a graph with $V(G) = \{v_1,\ldots, v_n\}$ and let $k_i$, $i\in [n]$, be positive integers. If $S$ is a gp-set of $G$ that induces a complete subgraph of $G\sr$, and  $\min\{k_i:\ v_i\in S\} \ge \max\{k_i:\ v_i\notin S\}$, then
$$\gp(G[K_{k_1},\ldots, K_{k_n}]) = \sum_{i:v_i\in S}k_i = \omega((G[K_{k_1},\ldots, K_{k_n]})\sr)\,.$$
\end{theorem}

\proof
Let $G$ and its gp-set $S$ be as stated in the theorem. Then $\gp(G) = \omega(G\sr)$  by Theorem~\ref{thm:lower-bound}. To simplify the notation, let $\widehat{G} = G[K_{k_1},\ldots, K_{k_n}]$ in the rest of the proof. Moreover, if a vertex $v_i\in V(G)$ expands to $R = K_{k_i}$ in $\widehat{G}$, and $\widehat{v}\in V(R)$, then we will write $v_i = g(\widehat{v})$. That is, if $\widehat{v}\in \widehat{G}$, then $g(\widehat{v})$ is the vertex of $G$ that expands to the complete subgraph of $\widehat{G}$ to which $\widehat{v}$ belongs.

If $\widehat{x}, \widehat{y}\in V(\widehat{G})$, $\widehat{x}\ne  \widehat{y}$, then by the construction of $\widehat{G}$ we infer that
\begin{equation}
\label{eq:distance-in-hat}
d_{\widehat{G}}(\widehat{x},\widehat{y}) =
  \begin{cases}
      1; & g(\widehat{x}) = g(\widehat{y})\,, \\
      d_{G}(g(\widehat{x}),g(\widehat{y})); & g(\widehat{x}) \ne g(\widehat{y})\,.
    \end{cases}
\end{equation}
By Theorem~\ref{thm:gpsets}, the components of $G[S]$ are complete subgraphs of $G$, denote them with $Q_1,\ldots, Q_r$. Then $\gp(G) = \sum_{i=1}^r |V(Q_i)|$. Since each vertex of $Q_i$ expands to a complete subgraph of $\widehat{G}$, the complete subgraph $Q_i$ expands to a complete subgraph of $\widehat{G}$, we will denote it with $\widehat{Q}_i$.

We first claim that $\widehat{S} = \bigcup_{i=1}^r V(\widehat{Q}_i)$ is a general position set of $\widehat{G}$. If $\widehat{x}\in V(\widehat{Q}_{i})$ and $\widehat{x}'\in V(\widehat{Q}_{i'})$, where $i, i'\in [r]$, $i\ne i'$, then $d_{\widehat{G}}(\widehat{x},\widehat{x}') = d_{G}(g(\widehat{x}),g(\widehat{x}'))$ holds by~\eqref{eq:distance-in-hat}. Therefore, since $\{Q_1,\ldots, Q_r\}$ form an in-transitive, distance-constant partition of $S$, the complete subgraphs $\{\widehat{Q}_1,\ldots, \widehat{Q}_r\}$ form an in-transitive, distance-constant partition of $\widehat{S}$. Hence, in view of Theorem~\ref{thm:gpsets}, $\widehat{S}$ is a general position set of $\widehat{G}$.

We next claim that $\widehat{S}$ is a gp-set of $\widehat{G}$. Assume on the contrary that there exists a general position set $\widehat{T}$ of $\widehat{G}$ such that $|\widehat{T}| > |\widehat{S}|$. Applying Theorem~\ref{thm:gpsets} again we know that the components of $\widehat{G}[\widehat{T}]$ are complete graphs. Let $T = \{g(\widehat{x}):\ \widehat{x}\in \widehat{T}\}$. Since $\widehat{T}$ is a general position set and because of~\eqref{eq:distance-in-hat} we infer that $T$ is a general position set of $G$. But since $\min\{k_i:\ v_i\in S\} \ge \max\{k_i:\ v_i\notin S\}$ and $|\widehat{T}| > |\widehat{S}|$ it follows that $|T| > |S| = \gp(G)$, a contradiction.

We have thus proved that $\gp(\widehat{G}) = \sum_{i:v_i\in S}k_i$. To complete the proof we need to show that also $\omega(\widehat{G}\sr) =  \sum_{i:v_i\in S}k_i$. Since $S$ is a complete subgraph of $G\sr$ and because of~\eqref{eq:distance-in-hat} we get that $\widehat{S}$ is a set of MMD vertices of $\widehat{G}$. By the equality part of Theorem~\ref{thm:lower-bound} we thus have $\omega(\widehat{G}\sr) = \gp(\widehat{G}) = \sum_{i:v_i\in S}k_i$.
\qed

\section{Rooted product graphs}

By a \emph{rooted graph} we mean a connected graph having one fixed vertex called the \emph{root} of the graph. Consider now a connected graph $G$ of order $n$, and let $H$ be a rooted graph with root $v$. The \emph{rooted product graph} $G\circ_v H$ is the graph obtained from $G$ and $n$ copies of $H$, say $H_1,\dots,H_n$, by identifying the root of $H_i$ with the $i^{\rm th}$ vertex of $G$, see~ \cite{godsil-1978, schwenk-1974}. To formulate the following result, the notion of an interval between vertices $u$ and $v$ of a graph $G$, defined as $I_G(u,v) = \{w:\ d_G(u,v) = d_G(u,w) + d_G(w,v)\}$, will be useful.

\begin{theorem}
\label{thm:rooted}
Let $G$ be any connected graph of order $n\ge 2$, and let $H$ be a rooted graph with root $v$.
\begin{itemize}
  \item[{\rm (i)}] $\gp(G\circ_v H)=n=\omega((G\circ_v H)\sr)$ if and only if $H$ is a path and $v$ is a leaf of $H$.
  \item[{\rm (ii)}] If $H$ contains a gp-set $S$ not containing $v$ and such that for each pair of vertices $u,w\in S$ neither $u\in I_H(v,w)$ nor $w\in I_H(v,u)$, then $\gp(G\circ_v H)=n\cdot \gp(H)$. Moreover, if in addition $S$ is a maximum clique in $H\sr$, then $\gp(G\circ_v H)=\omega((G\circ_v H)\sr)$.
  \item[{\rm (iii)}] Suppose $H$ is not a path rooted in one of its leaves. If every gp-set $S$ of $H$ either contains the root $v$, or contains two vertices $x,y$ such that ($x\in I_H(v,y)$ or $y\in I_H(v,x)$), then $2n\le \gp(G\circ_v H)\le n(\gp(H)-1)$. Particularly, if every gp-set of $H$ contains the root $v$, then $\gp(G\circ_v H)=n(\gp(H)-1)$.
\end{itemize}
\end{theorem}

\proof
(i) If $G$ is $P_2$ and $H$ is a path rooted in a leaf $v$, then $G\circ_v H$ is also a path, and so $\gp(G\circ_v H)=2=\omega((G\circ_v H)\sr)$. In this sense, from now we may assume $G$ is different from $P_2$.

If $H$ is a path and $v$ is a leaf, then clearly the set formed by the remaining leaves of all copies of $H$ forms a general position set of $G\circ_v H$, and so $\gp(G\circ_v H)\ge n$. Now, suppose $\gp(G\circ_v H)>n$ and let $S$ be gp-set of $G\circ_v H$. In consequence, by the pigeon hole principle there exists a copy, say $H_i$, of $H$ such that $|S\cap V(H_i)|\ge 2$, and indeed, it must happen $|S\cap V(H_i)|=2$. But then, the two vertices of $S\cap V(H_i)$ and any other distinct vertex of $S$ lie on a common geodesic, which is not possible (note that this third vertex always exists since $G$ is not $P_2$). Therefore, $\gp(G\circ_v H)\le n$ and the first equality follows. On the other hand, it can be easily observed that the strong resolving graph of $G\circ_v H$ is formed by a component isomorphic to a complete graph $K_n$,  and the remaining vertices of it are isolated ones. Thus, $\omega((G\circ_v H)\sr)=n$, which gives the second equality.

On the other hand, assume $\gp(G\circ_v H)=n=\omega((G\circ_v H)\sr)$. If $H$ is not a path rooted in one of its leaves, then there are at least two vertices of $H$, say $a,b$, such that $d_H(a,v)=d_H(b,v)$. In consequence, the set formed by the union of the copies of $a$ and $b$ in each copy of $H$ is a general position set of $G\circ_v H$ of cardinality $2n$, which is not possible. Thus, $H$ must be a path rooted in one of its leaves.

(ii) Let $A_i$, $i\in [n]$, be a gp-set of $H_i$ satisfying the statement of the item, and let $A=\bigcup_{i=1}^n A_i$. Then $A$ is a general position set of $G\circ_v H$, and so $\gp(G\circ_v H)\ge n\cdot \gp(H)$. Hence, suppose $\gp(G\circ_v H)> n\cdot \gp(H)$ and let $B$ be a gp-set of $G\circ_v H$. Thus, again by the pigeon hole principle, there must be a copy $H_j$ of $H$ such that $|B\cap V(H_j)|> \gp(H)$, but this is impossible since each copy of $H$ is an isometric subgraph of $G\circ_v H$ and $B\cap V(H_j)$ is a general position set of the graph induced by $H_j$. Consequently, $\gp(G\circ_v H)\le n\cdot \gp(H)$ and the equality follows.

On the other hand, assume that $A_i$ is a maximum clique in $H\sr$. Hence $\gp(H)=\omega(H\sr)$. Thus, from the above we get that $\gp(G\circ_v H)= n\cdot \omega(H\sr)$. It remains only to prove that $\omega((G\circ_v H)\sr)=n\cdot\omega(H\sr)$. Since any two vertices $u,w\in A_i$ satisfy that neither $u\in I_H(w,v)$ nor $w\in I_H(u,v)$, we see that $A=\bigcup_{i=1}^n A_i$ (defined as above) is also a clique in $(G\circ_v H)\sr$, and so $\omega((G\circ_v H)\sr)\ge n\cdot\omega(H\sr)$. Clearly, if we suppose that $\omega((G\circ_v H)\sr)> n\cdot\omega(H\sr)$, then we obtain that some $(H_j)\sr$ contains a clique of cardinality larger than $\omega(H\sr)$, which is not possible. Therefore, the required equality follows.

(iii) If every gp-set $S$ of $H$ either contains $v$ or contains two vertices $x,y$ such that without loss of generality $v$ belongs to an $x,y$-geodesic, then in order to construct a general position set of $G\circ_v H$ from the union of the gp-sets $S$ in each copy of $H$, we need to remove some vertices from each copy of $S$ including $v$ if it is the case. Clearly, the maximum number of vertices we may remove from $S$ is $|S|-1$, since a set formed by one vertex from each copy of $H$ is a general position set of $G\circ_v H$. However, as we next show by removing from $S$ at most $|S|-2$ vertices or removing $|S|-1$ and adding one other vertex not from $S$, we also obtain a general position set. Since $H$ is not a path rooted in one of it leaves, there are two vertices $x_i,y_i\in V(H_i)$ such that $d_{H_i}(x_i,v)=d_{H_i}(y_i,v)$, $i\in [n]$. Thus, the set $Q=\bigcup_{i=1}^n\{x_i,y_i\}$ is a general position set of cardinality $2n$ in $G\circ_v H$, and the lower bound follows. On the other hand, let $D$ be a gp-set of $G\circ_ v H$ and for every $i\in [n]$, let $D_i=D\cap V(H_i)$. If there is a set $D_j$ such that $|D_j|=\gp(H)$ (note that $|D_j|\le \gp(H)$ since $V(H_j)$ induces an isometric subgraph of $G\circ_v H$), then either $v\notin D_j$ or for any two vertices $x,y$ of $D_j$ it must happen that $v$ does not belong to an $x,y$-geodesic nor to a $y,x$-geodesic, but this is a contradiction with our assumption. Therefore, for every $i\in [n]$, $|D_i|\le \gp(H)-1$, which implies the upper bound.

We now consider the particular case in which every gp-set of $H$ contains the root $v$. Let $S_i$ be a gp-set of the copy $H_i$ of $H$ and let $S=\bigcup_{i=1}^n (S_i\setminus\{v\})$. Since $v$ belongs to $S_i$, it happens that $v$ does not belong to any $x,y$-geodesic for every $x,y\in S\setminus\{v\}$. Thus, no three vertices of $S$ lie on the same geodesic of $G\circ_v H$, and so, $S$ is a general position set of $G\circ_v H$. Therefore, $\gp(G\circ_v H)=n(\gp(H)-1)$.
\qed

Equality in the lower bound given in item (iii) of Theorem~\ref{thm:rooted} above can be noticed if $H$ is a path rooted in a vertex of degree two, where $\gp(G\circ_v H)=2n$. On the other hand, as we next observe the value of $\gp(G\circ_v H)$ can be very far from both bounds given above.

\begin{proposition}
There is a graph $G$ of order $n$ and a graph $H$ rooted in a vertex $v$ such that $n \ll \gp(G\circ_v H) \ll n(\gp(H)-1)$.
\end{proposition}

\proof
We consider a graph $H$ obtained as follows. We begin with a complete graph $K_r$. Next we add a vertex $v$ and join it by an edge to exactly $t$ vertices of $K_r$, where $2\le t\le r-2$, and choose $v$ as the root of this graph. Note that $H$ has only one gp-set $S$ formed by the set of vertices of the complete graph $K_r$. Also, note that if $x$ is adjacent to $v$, then $x\in I_H(y,v)$ for any $y$ not adjacent to $v$.

Now, let $G$ be a connected graph and let $A$ be a gp-set of $G\circ_v H$. Thus, if $A_i=A\cap V(H_i)$, then either every vertex of $A_i$ is adjacent to $v$ or no vertex of $A_i$ is adjacent to $v$, and so $|A_i|\le \max\{t,r-t\}$. As a consequence, $\gp(G\circ_v H)=|A|=\sum_{i=1}^{n}|A_i|\le n\cdot\max\{t,r-t\}$. On the other hand, the union of all neighbors of $v$ in each copy of $H$ in $G\circ_v H$, or the union of all not neighbors of $v$ in each copy of $H$ in $G\circ_v H$ is clearly a general position set of $G\circ_v H$, and so $\gp(G\circ_v H)\ge n\cdot\max\{t,r-t\}$, which implies the equality $\gp(G\circ_v H)=n\cdot\max\{t,r-t\}$. Since $\gp(H)=r$, the difference $n(\gp(H)-1)-\gp(G\circ_v H)=n(r-1-\max\{t,r-t\})$ can be arbitrarily large, as well as the difference $\gp(G\circ_v H)-n=n(\max\{t,r-t\}-1)$.
\qed

\section*{Acknowledgements}

We acknowledge the financial support from the Slovenian Research Agency (research core funding No.\ P1-0297 and projects J1-9109, N1-0095, N1-0108). This research was done while the second author was visiting the University of Ljubljana, Slovenia, supported by ``Ministerio de Educaci\'on, Cultura y Deporte'', Spain, under the ``Jos\'e Castillejo'' program for young researchers (reference number: CAS18/00030).


\end{document}